\documentclass[11pt,reqno]{article}
\usepackage[all]{xy}
\usepackage{amsmath}
\usepackage{amsthm}
\usepackage{amscd}
\usepackage{amssymb}

\input{epsf.tex}

\newcommand{\nb}[1]{#1\nobreakdash-}

\theoremstyle{definition}

\theoremstyle{plain}
\newtheorem{theorem}{Theorem}
\newtheorem{proposition}[theorem]{Proposition}

\newtheorem{corollary}[theorem]{Corollary}

\newtheorem{conjecture}[theorem]{Conjecture}
\newtheorem{fact}[theorem]{Fact}
\newtheorem*{theorem*}{Theorem}

\newcounter{remarks}
{\paragraph*{Remarks}\smallskip
     \begin{list}{\arabic{remarks}. }{\usecounter{remarks}%
          \setlength{\leftmargin}{0in}%
          \setlength{\rightmargin}{0in}%
          \setlength{\labelsep}{0pt}%
          \setlength{\labelwidth}{0pt}%
          \setlength{\listparindent}{0pt}%
     }
}
{
\end{list}
}

\newcommand\ds\displaystyle
\newcommand\wt[1]{\widetilde{#1}}

\DeclareMathOperator{\Out}{Out}

\DeclareMathOperator{\vcd}{vcd}
\DeclareMathOperator{\cd}{cd}

\DeclareMathOperator{\Length}{Length}

\DeclareMathOperator{\SL}{SL}
\DeclareMathOperator{\GL}{GL}

\DeclareMathOperator{\QIMap}{\widehat{\QI}}
\DeclareMathOperator{\QC}{QC}

\DeclareMathOperator{\Isom}{Isom}

\DeclareMathOperator{\Aut}{Aut}
\DeclareMathOperator{\Stab}{Stab}

\DeclareMathOperator\Conf{Conf}

\DeclareMathOperator\diam{diam}

\DeclareMathOperator\PD{PD}
\DeclareMathOperator\VN{VN}

\DeclareMathOperator\supp{supp}

\newcommand\R{{\mathbf R}}
\newcommand\reals{\R}

\renewcommand\H{{\mathbf H}}
\newcommand\Hyp\H
\newcommand\hyp{\H}

\newcommand\Z{{\mathbf Z}}

\newcommand\solv{{\ensuremath{\text{\scshape solv}}}}
\newcommand\Solv\solv

\newcommand\inject{\hookrightarrow}

\newcommand\infinity{\infty}
\newcommand\bndry{\partial}
\newcommand{\bdy}{\bndry}
\newcommand{\from}{\colon}
\def\composed{\circ}
\newcommand\suchthat{\bigm|}
\newcommand\inverse{{-1}}
\newcommand\inv{\inverse}

\newcommand\union{\cup}

\newcommand\absvalue[1]{{\left| #1 \right|}}
\newcommand\abs[1]{\absvalue{#1}}

\newcommand\Id{\text{Id}}
\newcommand\A{\mathcal A}
\newcommand\intersect{\cap}

\newcommand\subgroup{<}
\newcommand\semidirect{\rtimes}

\newcommand\Teichmuller{Teichm\"uller}
\newcommand\Poincare{Poincar\'e}

\newcommand\MCG{\mathcal{MCG}}

\DeclareMathOperator\QI{QI}
\newcommand\cross{\times}

\newcommand\Haus{{\mathcal H}}
\newcommand\G{{\mathcal G}}

\renewcommand\O{{\mathcal O}}

\newcommand\X{{\mathcal X}}

\newcommand\fol{{\mathcal F}}

\newcommand\F\fol

\newcommand\C{\mathcal C}
\newcommand\Teich{\mathcal T}
\newcommand\T\Teich

\DeclareMathOperator\Comm{Comm}

\newcommand\e\epsilon
\newcommand\<\langle
\renewcommand\>\rangle

\newcommand\ad{\text{ad}}

\newcommand\suf{{\text{suf}}}
\newcommand\uf{{\text{uf}}}

\title{Homology and dynamics in quasi-isometric rigidity\\ of once-punctured mapping
class groups}
\author{Lee Mosher}

\begin{document}
\maketitle

\begin{abstract}
In these lecture notes, we combine recent homological methods of Kevin Whyte with older
dynamical methods developed by Benson Farb and myself, to obtain a new quasi-isometric
rigidity theorem for the mapping class group $\MCG(S_g^1)$ of a once punctured surface
$S_g^1$: if $K$ is a finitely generated group quasi-isometric to $\MCG(S_g^1)$ then
there is a homomorphism $K \to \MCG(S_g^1)$ with finite kernel and finite index image.
This theorem is joint with Kevin Whyte.
\end{abstract}

Gromov proposed the program of classifying finitely generated groups according to their
large scale geometric behavior. The goal of this paper is to combine recent
homological methods of Kevin Whyte with older dynamical methods developed by Benson
Farb and myself, to obtain a new quasi-isometric rigidity theorem for mapping class
groups of once punctured surfaces:

\begin{theorem}[Mosher-Whyte] 
\label{TheoremMCGQIRigidity}
If $S^1_g$ is an oriented, once-punctured surface of genus
$g \ge 2$ with mapping class group $\MCG(S_g^1)$, and if $K$ is a finitely generated
group quasi-isometric to $\MCG(S_g^1)$, then there exists a homomorphism $K \to
\MCG(S_g^1)$ with finite kernel and finite index image.
\end{theorem}

This theorem will be restated later with a more quantitatively precise conclusion;
see Theorem~\ref{TheoremMCGQIQuant}. 

Whyte is also able to apply his techniques to obtain a strong quasi-isometric
rigidity theorem for the group $\Z^n\semidirect\GL(n,\Z)$, which we will not state
here.

Our theorem about $\MCG(S_g^1)$, answers a special case of:

\begin{conjecture} If $S$ is a nonexceptional surface of finite type then for any
finitely generated group $K$ quasi-isometric to $\MCG(S)$ there exists a homomorphism
$K\to\MCG(S)$ with finite kernel and finite index image.
\end{conjecture}

The exceptional surfaces that should be ruled out include several for which we already
have quasi-isometric rigidity theorems of a different type: the sphere with $\le 3$
punctures whose mapping class groups are finite; the once-punctured torus and the four
punctured sphere whose mapping class groups are commensurable to a free group of rank
$\ge 2$. Probably the techniques used for the once-punctured case may not be too
useful in the general case.

The theorem about $\MCG(S_g^1)$, and Whyte's results about $\Z^n \semidirect
\GL(n,\Z)$, are both about ``universal extension'' groups of certain $\PD(n)$ groups:
$\MCG(S_g^1) \approx \Aut(\pi_1 S_g)$ is the universal extension of the
$\PD(2)$ group $\pi_1 S_g$; and $\Z^n \semidirect \GL(n,\Z)$ is the universal extension
of the $\PD(n)$ group $\Z^n$. If one wishes to pursue quasi-isometric rigidity for the
group $\Aut(F_n)$, where $F_n$ is the free group of rank $\ge 2$, noting that
$\Aut(F_n)$ is the universal extension of $F_n$, the difficulty is that the
homological techniques we shall use do not apply: the \Poincare\ duality groups $\pi_1
S_g$ and $\Z^n$ each have a fundamental class in uniformly finite homology, which $F_n$
does not have.

\paragraph{Contents:} This paper is based on \LaTeX\ slides that were prepared for
lectures given at the LMS Durham Symposium on Geometry and Cohomology in Group Theory,
July 2003. Here is an outline of the paper, based approximately on my four lectures at
the conference:

\begin{enumerate}
\item Survey of results and techniques in quasi-isometric rigidity.
\item Whyte's techniques: uniformly finite homology applied to extension
groups.
\item Surface group extensions and Mess subgroups.
\item Dynamical techniques: extensions of surface groups by pseudo-Anosov
homeomorphisms.
\end{enumerate}

\paragraph{Acknowledgements.} Supported in part by NSF grant DMS 0103208. Thanks to
the organizers of the LMS Durham Symposium, July 2003.

\section{Survey of results and techniques in QI-rigidity}

A map $f \from X \to Y$ of metric spaces is a \emph{quasi-isometric embedding} if
$\exists K \ge 1, C \ge 0$ such that
$$\frac{1}{K}  \cdot d_X(x,y) - C \le d_Y(fx,fy) \le K \cdot d_X(x,y) + C
$$
A \emph{coarse inverse} for $f$ is a quasi-isometry $\bar f \from Y \to X$ s.t.
$$d_{\sup}(\bar f \composed f, \Id_X), \qquad d_{\sup}(f \composed \bar f,
\Id_Y) < \infinity
$$
A coarse inverse exists if and only if $\exists C' \ge 0$ such that $\forall y \in Y$
$\exists x \in X$ such that 
$$d_Y(fx,y) \le C'
$$
If this happens then $f \from X \to Y$ is a \emph{quasi-isometry},
and $X,Y$ are \emph{quasi-isometric} metric spaces. We will use the abbreviation ``QI''
for ``quasi-isometric''.

Given a finitely generated group $G$, a \emph{model space} for $G$ is a metric space
$X$ on which $G$ acts by isometries such that:
\begin{itemize}
\item $X$ is \emph{proper}, meaning that closed balls are compact.
\item $X$ is \emph{geodesic}, meaning that any $x,y \in X$ are connected by a
rectifiable path $\gamma$ such that $\Length(\gamma) = d(x,y)$.
\item The action is properly discontinuous and cobounded.
\end{itemize}

Examples of model spaces:
\begin{itemize}
\item The Cayley graph of $G$ with respect to a finite generating set.
\item $X = \wt Y$ where $Y$ is a compact, connected Riemannian manifold or piecewise
Riemannian cell complex, and $G = \pi_1 Y$. 
\end{itemize}

\paragraph{Fact:} If $X,Y$ are two model spaces for $G$ then $X,Y$ are
quasi-isometric. Also, any model space is quasi-isometric to $G$ with its word metric.

As a consequence, a finitely generated group $G$ has a notion of geometry that is
well-defined up to quasi-isometry, namely the geometry of any model space, or of $G$
itself with a word metric.

\paragraph{Definition:} Two finitely generated groups are \emph{quasi-isometric} if,
when equipped with their word metrics, they become quasi-isometric as metric spaces;
equivalently, their Cayley graphs are quasi-isometric.

\paragraph{Notation:} Given $\G$ a collection of finitely generated groups,  let 
$\<\G\> = \{$all groups quasi-isometric to some group in $\G\}$. More generally, given 
$\X =$ a collection of metric spaces, let $\<\X\> =\{$all groups quasi-isometric to
some metric space in $\X\}$.

\paragraph{Examples of QI-rigidity theorems.} To reformulate Gromov's program
in a practical way: given a collection of metric spaces $\X$, describe the collection
of groups $\<\X\>$, preferably in simple algebraic or geometric terms that do not
invoke the concept of quasi-isometry. Also, describe all of the quasi-isometry classes
within~$\<\X\>$. In particular, identify interesting classes of groups $\G$ that are
\emph{QI-rigid}, meaning $\G=\<\G\>$.  There are many theorems describing interesting
QI-rigid classes of groups, proved using an incredibly broad range of mathematical
tools. 

\paragraph{Example:} Gromov's polynomial growth theorem
\cite{Gromov:PolynomialGrowth} implies

\begin{theorem} The class of virtually nilpotent groups is quasi-isometrically rigid.
The class of virtually abelian groups is quasi-isometrically rigid, with one QI-class
for each rank.
\end{theorem}

Within the class of virtually nilpotent groups, there are many interesting
QI-invariants:
\begin{itemize}
\item  The Hirsch rank is a QI-invariant. 
\item The sequence of ranks of the abelian subquotients is a (finer) QI-invariant. 
\item There is an even finer QI-invariant of a virtually nilpotent Lie group $G$:
Pansu proved that the asymptotic cone of $G$ is a graded Lie group, whose associated
graded Lie algebra is a quasi-isometry invariant that subsumes the previous
invariants \cite{Pansu:croissance}. 
\item This is still not the end of the story: recently Yehuda Shalom produced two
finitely generated nilpotent groups which are not quasi-isometric but whose associated
graded Lie algebras are isomorphic \cite{Shalom:amenable}.
\end{itemize}
The full QI-classification of virtually nilpotent groups remains unknown.

\paragraph{Example:} Stallings' ends theorem \cite{Stallings:ends} implies

\begin{theorem} The class of groups which splits over a finite group is
quasi-isometrically rigid. For each $n \ge 2$, the class $\<F_n\>$ consists of all
groups that are virtually free of rank $\ge 2$.
\end{theorem}

By work of Papasoglu and Whyte \cite{PapasogluWhyte:ends}, combined with Dunwoody's
accessibility theorem \cite{Dunwoody:Accessible}, the QI classification of
finitely presented groups that split over a finite group is completely reduced to the
QI classification of one ended groups. 

\paragraph{Example:} Sullivan proved \cite{Sullivan:ErgodicAtInfinity} that any
uniformly quasiconformal action on $S^2$ is quasiconformally conjugate to a conformal
action. This implies:

\begin{theorem}[Sullivan--Gromov] $\<\hyp^3\>$ consists of all groups $H$ for which
there exists a homomorphism $H \to \Isom(\hyp^3)$ with finite kernel and whose image
is a cocompact lattice. 
\end{theorem}

This theorem is prototypical of a broad range of QI-rigidity theorems, including our
theorem about $\MCG(S_g^1)$. However, the conclusion of our theorem should be
contrasted with the Sullivan--Gromov theorem: the latter gives only a ``topological''
characterization of $\<\hyp^3\>$, which does not serve to give us an effective list of
those groups in $\<\hyp^3\>$. There is still no effective listing of the cocompact
lattices acting on $\hyp^3$. The conclusion of our theorem gives an ``algebraic''
characterization of $\<\MCG(S_g^1)\>$, allowing an effective listing.

\paragraph{Example:} Rich Schwartz proved a strong quasi-isometric rigidity theorem for
noncocompact lattices in $\Isom(\hyp^3)$ \cite{Schwartz:RankOne}. To state the theorem
we need some definitions.

\subparagraph{The commensurator group:} Given two groups $G,H$, a \emph{commensuration}
from $G$ to $H$ is an isomorphism from a finite index subgroup of $G$ to a finite
index subgroup of $H$. Two commensurations are \emph{equivalent} if they agree upon
restriction to another finite index subgroup. The \emph{commensurator group} $\Comm(G)$
is the set of self-commensurations of $G$ up to equivalence, with the following group
law: given commensurations $\phi \from A \to B, \psi \from C \to D$  restrict the
range of $\phi$ and the domain of $\psi$ to the finite index subgroup $B\intersect C$,
and then compose $\psi \composed \phi$. 

The left action of $G$ on itself by conjugation induces a homomorphism $G \to
\Comm(G)$, whose kernel is the \emph{virtual center} of $G$, consisting of all elements
$g \in G$ such that the centralizer of $g$ has finite index in $G$. 

Two groups $G,H$ are \emph{abstractly commensurable} if there exists a commensuration
from $G$ to $H$. Any abstract commensuration from $G$ to $H$ induces an isomorphism
from $\Comm(G)$ to $\Comm(H)$.

\begin{theorem}[Schwartz] If $G$ is a noncocompact, nonarithmetic lattice in
$\Isom(\hyp^3)$, then $\<G\>$ consists of those finitely generated groups $H$ which are
abstractly commensurable to $G$. More precisely, the homomorphism $G \to \Comm(G)$ is
an injection with finite index image, and $\<G\>$ consists of those finitely generated
groups $H$ for which there exists a homomorphism $H' \to \Comm(G)$ with finite kernel
and finite index image.
\end{theorem}

This theorem gives a very precise and effective enumeration of $\<\G\>$, similar to the
conclusion of our main theorems. Schwartz' theorem also can be formulated in the
arithmetic case, although there the homomorphism $G \to \Comm(G)$ has infinite index
image.

The general techniques of Sullivan-Gromov theorem, and of Schwartz' theorem give
models for the proof of our main theorem, as we now explain.

\paragraph{Technique: the quasi-isometry group of a group.} Consider a metric space
$X$, for example a model space for a finitely generated group. Let $\QIMap(X)$ be the
set of self quasi-isometries of $X$, equipped with the operation of composition.
Define an equivalence relation on $\QIMap(X)$, where $f \sim g$ if
$d_{\sup}(f,g) = \sup\{fx,gx\} < \infinity$. Composition
descends to a group operation on the set of equivalence classes, giving a group
$$\QI(X) = \text{the \emph{quasi-isometry group} of $X$}
$$
Notation: let $[f]$ denote the equivalence class of $f$ in $\QI(X)$. Note that
$[f]^\inv = [\bar f]$ for any coarse inverse $\bar f$ to $f$.

For any quasi-isometry $f \from X \to Y$ we obtain an isomorphism $\ad_f \from \QI(X)
\to \QI(Y)$ defined by $\ad_f[g] = [f \composed g \composed \bar f]$, where $\bar f$
is any coarse inverse for $f$. 

It follows that if $G$ is a finitely generated group then the \emph{quasi-isometry
group} of $G$ is well-defined up to isomorphism by taking it to be $\QI(X)$ for any
model space $X$ of $G$. 

The group $\QI(G)$ is an important quasi-isometry invariant of a group $G$, and it is
often important to be able to compute it. Here are some properties of $\QI(G)$,
followed a little later by some examples of computations.

The left action of $G$ on itself by multiplication, defined by $L_g(h) = gh$, induces a
homomorphism $G\to\QI(G)$ whose kernel is the virtual center. Also, the left action of
$G$ on itself by conjugation, defined by $C_g(h) = ghg^\inv$, also defines a
homomorphism $G\to\Comm(G)$. These two homomorphisms are identical, because
$d_{\sup}(L_g,C_g)$ equals the word length of $g$.

Every commensuration defines a natural quasi-isometry of $G$, well defined in
$\QI(G)$ up to equivalence of commensurations, thereby defining a homomorphism
$\Comm(G) \to \QI(G)$. The homomorphism $G \to \QI(G)$ factors as
$$G \to \Comm(G) \to \QI(G)
$$

\paragraph{Technique: quasi-actions} Let $G$ be a finitely generated group, $X$ a
model space for $G$, and $H$ a finitely generated group quasi-isometric to $G$. Fix a
quasi-isometry $\Phi \from H \to X$ and a coarse inverse $\bar \Phi \from X \to H$.
Define $A \from H \to \QIMap(X)$ by the formula 
$$A(h) = \Phi \composed L_h
\composed \bar\Phi
$$ 
This map has the following properties:
\begin{description}
\item[$A$ is a quasi-action:] There exists constants $K \ge 1$, $C \ge 0$ such that
\begin{itemize}
\item The maps $A(h)$ are $K,C$ quasi-isometries for all $h \in H$
\item $d_{\sup}(A(hh'), A(h) \composed A(h'))  \le C$ for all $h,h' \in H$
\item $d_{\sup}(A(\Id),\Id) \le C$
\end{itemize}
and so we obtain a homomorphism $A \from H \to \QI(X)$.
\item[$A$ is proper:]  $\forall r \ge 0$ $\exists n$ such that if
$B,B'
\subset X$ have diameter $\le r$ then $$\abs{\{h \in H \suchthat (A(h) \cdot B)
\intersect B'
\ne \emptyset\}} \le  n$$
\item[$A$ is cobounded:] $\exists s \ge 0$ such that $\forall x,y \in X$ $\exists h
\in H$ such that $d(A(h) \cdot x,y) \le s$.
\end{description}

Given a group $G$ and a model space $X$, a common strategy in investigating
quasi-isometric rigidity of $G$ is: compute $\QI(X)$; and then describe those
homomorphisms $H \to \QI(X)$ arising from quasi-actions, called ``uniform''
homomorphism. If necessary, restrict to proper, cobounded quasi-actions. Try to
``straighten'' any such quasi-action.

\paragraph{Examples of QI-rigidity.} Here are some examples of how this strategy is
carried out, taken from the above examples:

\begin{proof}[Proof sketch for the Sullivan--Gromov rigidity theorem] For groups in the
quasi-isometry class of $\hyp^3$, the boundary is $\bdy\hyp^3=S^2$. 

First, one calculates $\QI(\hyp^3) = \QC(S^2)$, the group of quasi-conformal
homeomorphisms of $S^2$; this is a classical result in quasiconformal geometry.

Second, the isometry group $\Isom(\hyp^3) = \Conf(S^2)$ is a uniform subgroup of
$\QC(S^2)$, and one proves that every uniform subgroup can be conjugated into
$\Conf(S^2)$. In other words, every quasi-action on $\hyp^3$ is quasiconjugate to an
action. This result, due to Sullivan, is the heart of the proof.

The properties of ``properness'' and ``coboundedness'' are invariant under
quasiconjugacy. It follows that if $H$ is a finitely generated group quasi-isometric
to $\hyp^3$ then $H$ has a proper, cobounded action on $\hyp^3$. In other words, there
is a homomorphism $H \to \Isom(\hyp^3)$ with finite kernel and discrete, cocompact
image. 
\end{proof}

By constrast we now give:

\begin{proof}[Proof sketch for the Schwartz rigidity theorem] Let $G$ be a noncocompact
lattice in $\hyp^3$. 

The heart of the proof is essentially a calculation
$$\QI(G) \approx \Comm(G)
$$
This calculation holds in both the arithmetic case and the nonarithmetic case, the
difference being that the induced map $G \to \Comm(G)$ has finite index image if and
only if $G$ is nonarithmetic. Assuming this to be the case, it follows that the
homomorphism $\Comm(G) \to \QI(G)$ is an isomorphism and that the map $G \to \Comm(G)
\approx \QI(G)$ is an injection with finite index image.

Schwartz' proof is actually a bit more quantitative, as follows. If $G$ is
nonarithmetic then there exists an embedding $$\Comm(G) \inject \Isom(\hyp^3)$$
whose image $\Gamma$ is a noncocompact lattice containing $G$ with finite
index, so that the injection $G \to \Comm(G)$ agrees with the inclusion $G \inject
\Gamma$. The hard part of Schwartz' proof is to show the following:
\begin{itemize}
\item $\forall K \ge 1$, $C \ge 0$ $\exists A \ge 0$ such that if $\Phi \from G \to G$
is a $K,C$ quasi-isometry then there exists $\gamma \in \Gamma$ such that
$$d_{\sup}(\Phi,L_\gamma) \le A
$$
\end{itemize}
To be more precise, the sup distance on the left is a comparison of two different
functions from $G$ into $\Gamma$, one being $G \xrightarrow{\Phi} G \inject \Gamma$,
and the other being $G\inject \Gamma \xrightarrow{L_\gamma} \Gamma$.

Noting that any quasi-isometry of $G$ extends to the finite index supergroup $\Gamma$,
and that $\QI(G) \approx \QI(\Gamma)$, we can abstract this discussion as follows.

Consider a finitely generated group $\Gamma$, and suppose that the following holds:

\paragraph{Strong QI-rigidity:} $\forall K \ge 1$, $C \ge 0$ $\exists A \ge 0$ such
that if $\Phi \from \Gamma
\to \Gamma$ is a $K,C$ quasi-isometry then there exists $\gamma \in \Gamma$ such that
$$d_{\sup}(\Phi,L_\gamma) \le A
$$

This property, coupled with triviality of the virtual center (true for lattices
in $\Isom(\hyp^3)$ as well as for $\MCG(S_g^1)$), immediately imply that the
homomorphism $\Gamma\to\QI(\Gamma)$ is an isomorphism.

To complete the proof of Schwartz' Theorem, we now apply the following fact:

\begin{proposition} If $\Gamma$ is a strongly QI-rigid group whose virtual center is
trivial, then for any finitely generated group $H$ quasi-isometric to $\Gamma$ there
exists a homomorphism $H \to\Gamma$ with finite kernel and finite index image.
\end{proposition}

\begin{proof} As explained earlier, the left action of $H$ on itself
by translation can be quasiconjugated to a proper, cobounded quasi-action of $H$ on
$\Gamma$, which induces a homomorphism $\phi\from H\to \QI(\Gamma) = \Gamma$. 

Let $K \ge 1$, $C\ge 0$
be uniform constants for the quasi-action of $H$ on $\Gamma$. 

Applying strong QI-rigidity of $\Gamma$, we obtain a constant $A$ such that the
(quasi-)action of each $h\in H$ on $\Gamma$ is within sup distance $A$ of left
multiplication by $\phi(h)$. It immediately follows that the kernel of $\phi$ is
finite, because the quasi-action of $H$ is proper and so there are only finitely many
elements $h \in H$ for which $\phi(h)$ is within distance $A$ of the identity on
$\Gamma$.

It also follows that the image of $\phi$ has finite index, because the quasi-action of
$H$ on $\Gamma$ is cobounded, whereas the left action on $\Gamma$ of any infinite index
subgroup of $\Gamma$ is not cobounded.
\end{proof}

This completes the proof of Schwartz' Theorem.
\end{proof}

This proof immediately yields an interesting corollary:

\begin{corollary} If $\Gamma$ is strongly QI-rigid with trivial virtual center, then
every commensuration of $\Gamma$ is the restriction of an inner automorphism of
$\Gamma$.
\end{corollary}

Now I can give the more quantitative statement of the main theorem about $\MCG(S_g^1)$:

\begin{theorem}[Mosher-Whyte] 
\label{TheoremMCGQIQuant}
The group $\MCG(S_g^1)$ is strongly QI-rigid: for all $K \ge
1$, $C \ge 0$ there exists $A \ge 0$ such that for any $K,C$ quasi-isometry $\Phi
\from\MCG(S_g^1)\to\MCG(S_g^1)$ there exists $\gamma \in \MCG(S_g^1)$ for which
$d_{\sup}(\Phi,L_\gamma) < A$.
\end{theorem}

As an application, we get a new proof of a result of Ivanov:

\begin{corollary}[Ivanov] The injection $\MCG(S_g^1) \to \Comm(S_g^1)$ is an
isomorphism, that is, every commensuration of $\MCG(S_g^1)$ is the restriction of an
inner automorphism.
\end{corollary}

\section{Fiber preserving quasi-isometries}

In this section we explain Kevin Whyte's methods for using uniformly finite homology
classes and their supports to investigate quasi-isometric rigidity problems for
certain fiber bundles. Bruce Kleiner also outlined, at the AMS Ann Arbor conference
in 2000, how to use support sets to study quasi-isometric rigidity, using a coarse
version of the K\"unneth formula applied to fiber bundles.

Suppose one wants to investigate quasi-isometric rigidity for the fundamental group of
a graph of groups where each vertex and edge group is (virtually) $\pi_1$ of an
aspherical $n$-manifold for fixed integer $n \ge 0$. Let us focus on the example
$F_2\cross\Z^n$. 

To start the proof, pick a nice model space for $F_2 \cross \Z^n$, namely, $T \cross
\reals^n$ where $T$ is a Cayley tree of $F_2$. As described earlier, the technique
will be to study quasi-actions on $T \cross \reals^n$. The first step, carried out by
Farb and myself, is to prove that each quasi-isometry of $T \cross
\reals^n$ coarsely preserves the $\reals^n$ fibers:

\begin{theorem}[\cite{FarbMosher:ABC}]
\label{TheoremFM}
$\forall K,C$ $\exists A$ such that if $\Phi \from T \cross \reals^n \to T \cross
\reals^n$ is a $K,C$ quasi-isometry, then $\forall t \in T$ $\exists t' \in T$ such
that $d_\Haus(\Phi(t \cross \reals^n), t' \cross \reals^n) < A$.
\end{theorem}

The notation $d_\Haus(\cdot,\cdot)$ means Hausdorff distance between subsets of a
metric space.

More generally, this theorem is true when the product $T \cross \reals^n$ is replaced
by a ``coarse fibration'' over a tree whose fiber is a uniformly contractible
manifold, or even more generally by the Bass-Serre tree of spaces that arise from
a finite graph of groups whose vertex and edge groups are coarse $\PD(n)$ groups of
fixed dimension $n$ \cite{MSW:QTOne}.  

In this section we shall give Kevin Whyte's proof of Theorem~\ref{TheoremFM}. This
proof also applies to situations where the base space $T$ of the fibration is replaced
by certain higher dimensional complexes, for example: 
\begin{itemize}
\item Thick buildings.
\item The model space for $\MCG(S_g)$, over which a model space for $\MCG(S_g^1)$
fibers, with fiber $\hyp^2$. 
\item A model space for $\SL(n,\Z)$, over which a model space for $\Z^n \semidirect
\GL(n,\Z)$ fibers, with fiber $\reals^n$.
\end{itemize}

For the example $T \cross \reals^n$, the idea of the proof is that a subset of the
form $\text{(line in $T$)} \cross \reals^n \approx \reals^{n+1}$ is the support of a
``top dimensional uniformly finite homology class''. A quasi-isometry of $T \cross
\reals^n$ acts on such classes, coarsely preserving their supports. Each fiber is the
intersection of some finite number of these supports, and so the fibers are preserved.

In this proof it is necessary that the fiber be a uniformly contractible
\emph{manifold}, on which there is a ``uniformly finite'' fundamental class of full
support. (For those who live in outer space, that's why the proof does not apply to
the extension $1 \to F_n \to \Aut(F_n) \to \Out(F_n) \to 1$, whose fiber $F_n$ is not
a manifold and does not have a uniformly finite fundamental class of full support).

In order to make this proof rigorous, we have to discuss: 
\begin{enumerate}
\item Uniformly finite homology 
\item Top dimensional supports 
\item Application to fiber bundles
\end{enumerate}

\subsection{Uniformly finite homology}

Let $X$ be a simplicial complex. Fix a geodesic metric in which each simplex is a
regular Euclidean simplex with side length $1$. We say that $X$ is \emph{uniformly
locally finite} or ULF if $\exists A \ge 0$ such that the link of each simplex contains at
most $A$ simplices. We say that $X$ is \emph{uniformly contractible} or UC if $\forall
r>0$ $\exists s(r) \ge 0$ s.t. each subset $A \subset X$ with $\diam(A) \le r$ is
contractible to a point inside $N_s(A)$. The function $s(r)$ is called a \emph{gauge}
of uniform contractibility.

For example, any tree is uniformly contractible. Also, if $X$ is a contractible
simplicial complex, and if there is a cocompact, simplicial group action on $X$ (for
example if $X$ is the universal cover of a compact, aspherical simplicial complex),
the $X$ is UC and ULF.

\paragraph{Simplicial uniformly finite homology.} If $X$ is ULF, define $H^\suf_n(X)$:
a chain in $C_n^\suf(X)$ is a uniformly bounded assignments of integers to
$n$-simplices. Since $X$ is ULF, the boundary map $\bdy \from C_n^\suf(X) \to
C_{n-1}^\suf(X)$ is defined, and clearly $\bdy\bdy=0$.

\begin{theorem} 
\label{TheoremUF}
$H^\suf_*(X)$ is a quasi-isometry invariant among UC, ULF simplicial complexes. 
\end{theorem}

\begin{proof} We prove the theorem by defining \emph{uniformly finite homology}
$H^\uf_*(X)$ which is a large scale versioon of simplicial uniformly finite homology
$H^\suf_*(X)$, proving that $H^\uf_*(X)$ and $H^\suf_*(X)$ are isomorphic for UC, ULF
simplicial complexes, and then proving that $H^\uf_*(X)$ is a QI-invariant.

\paragraph{Step 1: Uniformly finite homology.} Define the $n$th Rips complex $R^n(X)$,
with one $k$-simplex for each ordered $k+1$-tuple of vertices with diameter $\le n$.
Note that $R^1(X)=X$. 

Since $X$ is ULF, it follows that $R^n(X)$ is ULF, and the sequence of inclusions
$$X=R^1(X) \subset R^2(X) \subset \cdots
$$
therefore induces homomorphisms
$$H^\suf_*(R^1(X)) \to H^\suf_*(R^2(X)) \to H^\suf_*(R^3(X)) \cdots
$$
Define the uniformly finite homology to be the direct limit
$$H^\uf_*(X) = \lim_{k \to \infinity} H^\suf_*(R^k(X))
$$
We can also describe $H^\uf_*(X)$ as the homology of a chain complex. We have a direct
system
$$C^\suf_*(R^1(X)) \to C^\suf_*(R^2(X)) \to C^\suf_*(R^3(X)) \to \cdots
$$
so we can take the direct limit 
$$C^\uf_*(X) = \lim_{k \to \infinity} C^\suf_*(R^k(X))
$$
The boundary homomorphism is defined, and the homology of this chain complex is
canonically isomorphic to $H^\uf_*(X)$. 

\paragraph{Step 2:} $H^\uf_*(X) = H^\suf_*(X)$. 

The identity map $i \from X \to X=R^1(X)$ induces a chain map $i \from C^\suf_*(X) \to
C^\uf_*(X)$. Using that $X$ is UC, we'll define a
chain map
$$j \from C^\uf_*(X) \to C^\suf_*(X)
$$
which is a uniform chain homotopy inverse to the inclusion $C^\suf_*(X) \to
C^\uf_*(X)$. 

The idea for defining $j$ is: connect the dots.

Consider a 1-simplex in $C^\uf_1(X)$, which means a 1-simplex $\sigma$ in $R^k(X)$ for
some $k$, which means $\sigma=(u,v)$ with $d(u,v) \le k$. Connecting the dots, we get
a 1-chain $j(\sigma)$ in $X$ with boundary $v-u$, consisting of at most $k$
1-simplices. Given now a simplicial uniformly finite 1-chain $\sum
a_\sigma \sigma$ in $R^k(X)$, the infinite sum 
$$j(\sum a_\sigma \sigma) = \sum a_\sigma j(\sigma)
$$ 
is defined because it is locally finite: for each simplex $\tau$ in $X$
there are finitely many terms of the sum $\sum a_\sigma j(\sigma)$ which assign a
nonzero coefficient to $\tau$. This finishes the definition of $j \from C^\uf_1(X) \to
C^\suf_1(X)$.

Next consider a 2-simplex $\sigma$ in $C^\uf_2(X)$, which means $\sigma = (u,v,w)$,
where $u,v,w$ have pairwise distances at most~$k$. Now connect the 2-dimensional dots:
the 1-chain 
$$j(u,v) + j(v,w) + j(w,u)
$$ 
is a cycle. Its support is a subset of diameter at most $3k/2$, and so
$$j(u,v) + j(v,w) + j(w,u) = \bdy j(\sigma)
$$ 
for some 2-chain $j(\sigma)$ supported on a subset of diameter at most $s(3k/2)$,
where $s$ is a gauge of uniform contractibility. More generally, the boundary of a
simplicially uniformly finite 2-chain in $R^2(X)$, is again defined as a locally
finite infinite sum. This finishes the definition of $j \from C^\uf_2(X) \to
C^\suf_2(X)$, and the chain map condition is obvious.

Now continue the definition of the chain map $j$ by induction, using connect-the-dots.

Similarly, using connect the dots and induction, we can construct a chain homotopy
between identity and $ji$, and similarly for $ij$. This finishes Step 2.

\paragraph{Step 3: QI-invariance of uniformly finite homology.} Consider a
quasi-isometry $\Phi\from X\to Y$ with coarse inverse $\bar \Phi \from Y \to X$, both
$K,C$ quasi-isometries, and $C$-coarse inverses of each other. Moving a bounded
distance, may assume $\Phi,\bar\Phi$ take vertices to vertices. 

If $d(u,v)=1$ then $d(\Phi u, \Phi v) \le p=K+C$. We therefore obtain an induced
simplicial map $X=R^1(X)\to R^p(X)$, inducing a chain map 
$$\Phi_\# \from C^\suf_*(R^1(X)) \to C^\suf_*(R^p(X))
$$ 
In the backwards direction, if $d(u',v') \le p$ in $Y$, then $d(\bar\Phi u', \bar\Phi
v') \le p' = Kp+C$, and so we get an induced simplicial map $R^p(Y)\to R^{p'}(X)$
inducing a chain map
$$\bar\Phi_\# \from C^\suf_*(R^{p}(Y)) \to C^\suf_*(R^{p'}(X))
$$ 
The composition $\bar\Phi \composed \Phi$ induces a chain map
$$\bar \Phi_\# \composed \bar \Phi_\# \from C^\suf_*(R^1(X)) \to
C^\suf_*(R^{p'}(X))
$$
but $\bar\Phi \composed \Phi$ is $C$-close to the identity map on vertices of $X$, and
so a connect-the-dots argument shows that this chain map is chain homotopic to the
inclusion map.

A similar argument applies to the composition $\Phi \composed \bar\Phi$. 

This finishes the proof that $H^\suf_*$ is a QI invariant.
\end{proof}

\subsection{Top dimensional supports.} Suppose now that $X$ is a UC, ULF simplicial
complex of dimension $d$. There are no simplices of dimension $d+1$, and so each class
$c \in H^\uf_d(X)$ is represented by a unique $d$-cycle in $C^\suf_d(X)$, also denoted
$c$. Its support $\supp(c)$ is therefore a well-defined subset of~$X$.

\begin{proposition}
\label{PropClassesPreserved}
With $X$ as above, every quasi-isometry of $X$ coarsely respects supports of classes in
$H^\uf_d(X)$. More precisely: $\forall K,C$ $\exists A$ such that if $\Phi \from X \to
X$ is a $K,C$ quasi-isometry, and if $c \in H^\uf_d(X)$, then 
$$d_\Haus\bigl(\Phi(\supp(c)), \supp(\Phi_*(c))\bigr) < A
$$
\end{proposition}

Most QI-rigidity theorems have a similar step: find some collection $\C$ of objects in
the model space which are coarsely respected by quasi-isometries: $\forall K,C$
$\exists A$ such that for each $K,C$ quasi-isometry $\Phi \from X \to X$, and for each
object $c \in \C$ there exists an object $c' \in \C$ such that
$$d_\Haus(\Phi(c),c') < A
$$

\begin{proof} Moving $\Phi$ a bounded distance, we may assume that $\Phi$ takes
vertices to vertices. We get an induced chain map
$$\Phi_\# \from C^\suf(X) \to C^\suf(R^p(X))
$$
Note first that for all $c \in C^\suf(X)$, the subset $\supp(\Phi_\#(c))$ is contained
in a uniformly bounded neighborhood of $\Phi_\#(\supp(c))$, where the support of a
chain in $C^\suf(R^p(X))$ is simply the set of vertices occurring among the
summands in the chain. In other words, $\Phi_\#$ induces coarse inclusion of supports.

Compose with the connect-the-dots map
$$C^\suf(R^p(X)) \to C^\suf(X)
$$
which also induces coarse inclusion of supports. By composition we obtain an induced
map
$$\Phi_{\#\#} \from C^\suf(X) \to C^\suf(X)
$$
which also induces coarse inclusion of supports. Since top dimensional supports are
unique, it follows that
$$\Phi(\supp(c)) \subset N_A(\supp(\Phi_{\#\#}(c)) = N_A(\supp(\Phi_*(c)))
$$
for some uniform constant $A$.

To get the inverse inclusion, applying the same argument to a coarse inverse $\bar
\Phi$ we have
$$\bar \Phi(\supp(\Phi_*(c))) \subset N_A(\supp(\bar \Phi_* \Phi_*(c))) = N_A(\supp(c))
$$
where the last equation follows from uniqueness of supports. Now apply $\Phi$ to both
sides of this equation:
\begin{align*}
\supp(\Phi_*(c)) &\subset N_{A'}(\Phi \bar \Phi(\supp(\Phi_*(c))) \\
 &\subset N_A'(\Phi(N_A(\supp(c))) \\
 &\subset N_{A''}(\Phi(\supp(c)))
\end{align*}
\end{proof}

\subsection{Application to fiber bundles.} Consider now a fiber bundle $\pi
\from E\to B$ with fiber $F_x$ over each $x \in B$. We assume that $E,B$ are UC, ULF
simplicial complexes, $\pi$ is a simplicial map, each fiber $F_x = \pi^\inv(x)$ is a
manifold of dimension $n$, and for each vertex $x$ the subcomplex $F_x$ is UC, with
gauge independent of~$x$. It follows that for each $k$-simplex $\sigma$, the $k+n$
simplices of $\pi^\inv(\sigma)$ that are not contained in $\pi^\inv(\bdy\sigma)$,
intersected with the fiber $F_\sigma = \pi^\inv(\text{barycenter}(\sigma))$, define a
cellular structure on $F_\sigma$ which is UC, with gauge independent of~$\sigma$.

Let $d = \dim(B)$, $n = \dim(F)$, $d+n = \dim(E)$.

Make the following assumption about the top
dimensional, uniformly finite homology $H^\uf_d(B)$: 

\begin{description}
\item[Top dimensional classes in $B$ separate points:] 
$\exists r>0$ so that $\forall s>0$ $\exists  D>0$ so that for any $x, y \in B$
with $d(x,y)>D$, there is a top dimensional class $c \in H^\uf_d(B)$ such that 
$$d(\supp(c),x) \le r \quad\text{and}\quad d(\supp(c),y) > s
$$
\end{description}

Example: $T\cross\reals^n$, where $T =$ Cayley tree of $F_2$. In the base space $T$,
each bi-infinite line is the support of a top dimensional class, and lines in $T$
clearly separate points.

Example to come: later we shall construct a model space for $\MCG(S_g^1)$, fibering
over model space for $\MCG(S_g)$, with fiber $\hyp^2$. The dimension of the base space
will equal $4g-5$, which is the virtual cohomological dimension of $\MCG(S_g)$. We
shall prove that the top dimensional classes of uniformly finite homology separate
points of $\MCG(S_g)$.

The main result for this section is that every quasi-isometry of the total space $E$
coarsely preserves fibers:

\begin{theorem}[Whyte] 
\label{TheoremFibers}
Consider the fibration $F \to E \to B$ as above, and assume that
top dimensional classes in $H^\uf_d(B)$ coarsely separate points. For all $K,C$ there
exists $A$ such that if $\Phi \from E \to E$ is a $K,C$ quasi-isometry, then
for each $x \in B$ there exists $x' \in B'$ such that
$$d_\Haus(\Phi(F_x), F_{x'}) \le A
$$
\end{theorem}

\begin{proof} The key observation of the proof is that the support of every top dimensional class
in $E$ is saturated by fibers. To be precise: for every top dimensional class
$c \in H^\uf_{d+n}(E)$, there exists a unique top dimensional class $c'=\pi(c)\in
H^\uf_d(B)$, such that
$$\supp(c) = \pi^\inv(\supp(c'))
$$

To see why, for each $d$-simplex $\sigma \subset B$, $\pi^\inv(B)$ is a
manifold with boundary of dimension $d+n$. So, for any class $c$ of dimension $d+n$,
if $\supp(c)$ contains some $d+n$ simplex in $\pi^\inv(\sigma)$, it follows that
$\supp(c)$ contains all of $\pi^\inv(\sigma)$.

The converse is also true: for each top dimensional class $c' \in H^\uf_d(B)$ there
exists a top dimensional class in $E$, denoted $c = \pi^\inv(c') \in H^\uf_{d+n}(E)$
such that $\supp(c) = \pi^\inv(\supp(c'))$: over each simplex $\sigma \subset
\supp(c')$, weight all the simplices in $\pi^\inv(\sigma)$ with the same weight as
$\sigma$, using a coherent orientation of fibers to choose the sign.

Thus, the projection $\pi$ induces an isomorphism 
$$H^\uf_{d+n}(E) \to H^\uf_d(B)
$$
so that a $d+n$-cycle $c$ in $E$, and the corresponding $d$-cycle $c'$ in $B$, are
related by
$$\supp(c) = \pi^\inv(\supp(c'))
$$

Now we use the property that supports of top dimensional classes in $B$ separate
points. Up to changing constants, it follows that:

\begin{description}
\item[Supports of top dimensional classes in $E$ separate fibers:] there exists $r>0$
such that for all $s>0$ there exists $D>0$ such that given fibers $F_x,F_y$ with
$d_\Haus(F_x,F_y) > D$, there is a top dimensional class $c \in H^\uf_{d+n}(E)$ so
that $F_x \subset N_r(\supp(c))$ but $F_y \intersect N_r(\supp(c)) = \emptyset$.
\end{description} 

From this it follows that fibers in $E$ are coarsely respected by a quasi-isometry.
Here are the details.

Fix a $K,C$ quasi-isometry $\Phi \from E \to T$ and a fiber $F_x$, $x \in G$. We want
to show that $\Phi(F_x)$ is uniformly Hausdorff close to some fiber $F_{x'}$. 

Fix some $R>r$, to be chosen later, and let $\C_x$ denote the collection of classes $c
\in H^\uf_{d+n}(E)$ such that $F_x\subset N_R(\supp(c))$. From the fact that top
dimensional classes in $E$ separate fibers, it follows that $F_x$ has (uniformly)
finite Hausdorff distance from the set of points $\xi \in E$ such that $\xi \in
N_R(\supp(c))$ for all
$c \in \C_x$.

Notation: let $\hat \C_x = \{\Phi_\#(c) \suchthat c \in \C_x\}$, and let $\hat c =
\Phi_\#(c)$.

By applying Proposition~\ref{PropClassesPreserved}, the Hausdorff distance between
$\Phi(\supp(c))$ and $\supp(\hat c)$ is at most a constant $A$, for any $c \in \C_x$.
Thus for any $\hat c \in \hat\C_x$ we have
$$\Phi(F_x) \subset N_{R'}(\supp(\hat c))
$$
where $R' = KR+C+A$. Now we say how large to choose $R$, namely, so that $R'>r$.

It now follows that $\Phi(F_x)$ has (uniformly) finite Hausdorff distance from the set
$\F$ of points $\eta \in E$ such that $\eta \in N_{R'}(\supp(\hat c))$ for all $\hat c
\in\hat\C_x$. But the set $\F$ is clearly a union of fibers of $E$. Pick one fiber
$F_{x'}$ in $\F$. Taking $s=R'$ in the definition of coarse separation of fibers,
there is a resulting $D$. If $F_{y'}$ is a fiber whose distance from $F_{x'}$ is more
than $D$, it follows that $F_{y'}$ is not contained in $\F$, because that would violate
coarse separation of fibers. This shows that the set $\F$ contains $F_{x'}$ and is
contained in the $D$-neighborhood of $F_{x'}$, that is, $\F$ has Hausdorff distance at
most $D$ from $F_{x'}$. But $\F$ also has (uniformly) finite Hausdorff distance from
$\Phi(F_x)$, and so $\Phi(F_x)$ has (uniformly) finite Hausdorff distance from
$F_{x'}$.

This finishes the proof of Theorem~\ref{TheoremFibers}.
\end{proof}

\section{Surface group extensions and Mess subgroups}

Let $S_g$ be a closed, oriented surface of genus $g \ge 2$, and let $S_g^1$ be $S_g$
minus a single base point $p$. There is a short exact sequence
$$1 \to \pi_1(S_g) \to \MCG(S_g^1) \to \MCG(S_g) \to 1
$$
The homomorphism $\MCG(S_g^1) \to \MCG(S_g)$ is the map that ``fills in the
puncture''. The homomorphism $\pi_1(S_g) \to \MCG(S_g^1)$ is the ``push'' map, which
isotopes the base point $p$ around a loop, at the end of the isotopy defining a map of
$S_g$ taking $p$ to itself; then remove the base point to define a mapping class on
$S_g^1$.

The main theorem of this section is:

\begin{theorem}
\label{TheoremCosets}
Every quasi-isometry of $\MCG(S_g^1)$ coarsely preserves the system of cosets of
$\pi_1(S_g)$. 
\end{theorem}

The meaning of this theorem is that for each $K \ge 1$, $C \ge 0$ there exists $A \ge
0$ such that if $\Phi \from \MCG(S_g^1) \to \MCG(S_g^1)$ is a $K,C$ quasi-isometry,
then for each coset $C$ of $\pi_1(S_g)$ there exists a coset $C'$ such that
$d_\Haus(\Phi(C),C') \le A$.

The setup of Theorem~\ref{TheoremCosets} is to represent the short exact sequence by a
fibration 
$$\hyp^2 \to E \to B
$$
as above, where $E$ is a model space for $\MCG(S_g^1)$, and $B$ is a model space for
$\MCG(S_g)$. The theorem can then be translated into geometric terms by saying that
every quasi-isometry of $E$ coarsely preserves the fibers. 

Whyte's idea for proving Theorem~\ref{TheoremCosets} is to apply 
Theorem~\ref{TheoremFibers}, by using Mess subgroups of $\MCG(S_g)$ to provide top
dimensional classes in the uniformly finite homology of $B$ that coarsely separate
points.

\subsection{Dimension of $\MCG(S_g)$.} Given a contractible model space $B$ for a group
$G$, one would expect that the top dimension in which $H^\uf_n(B)$ is nontrivial would
be
$n=\vcd(G)$. So we need the following formula of John Harer:

\begin{theorem}[\cite{Harer:vcd}]
$$\vcd(\MCG(S_g)) = 4g-5
$$ 
\end{theorem}

\begin{proof} We will sketch Harer's original proof of the upper bound $\vcd(\MCG(S_g))
\le 4g-5$, and then we shall give Geoff Mess' proof of the lower bound
$\vcd(\MCG(S_g)) \le 4g-5$. Mess' proof will provide the basic ingredients we need to
investigate top dimensional uniformly finite homology classes in a model space for
$\MCG(S_g)$.

To prove $\vcd(\MCG(S_g)) \le 4g-5$, by using the short exact sequence, in which 
$$\vcd(\text{kernel}) = \dim(\hyp^2) = 2
$$
it suffices to prove
$$\vcd(\MCG(S_g^1)) \le 4g-3
$$
Harer constructs a contractible complex $K$ of dimension $4g-3$ which is a model space
for $\MCG(S_g^1)$: $K$ is the complex of ``filling arc systems'' of the once punctured
surface. Fixing a base point $p \in S_g$, a filling arc system is a system of
arcs $\A = \{A_i\}$ whose interiors are pairwise disjoint, whose ends are all located
at the base point $p$, so that for each component $C$ of $S_g - \union \{A_i\}$,
regarding $C$ as the interior of a polygon whose sides are arcs of $\{A_i\}$, the
number of sides of $C$ is at least $3$. Each filling arc system $\A$ can be refined by
adding more arcs until it is triangulated, and the number of such arcs is called the
\emph{defect} of $\A$. The complex $K$ that Harer uses has one cell $C_\A$ of dimension
$d$ for each isotopy class of filling arc systems $\A$ of defect $d$, and the boundary
of $C_\A$ consists of all cells $C_{\A'}$ such that $\A\subset\A'$. Harer used Strebel
differentials to prove contradictiblity of $K$, but a purely combinatorial proof was
given by Hatcher \cite{Hatcher:triangulations}.

\subparagraph{Remarks:} Harer does \emph{not} directly construct a $4g-3$ dimensional
model space for $\MCG(S_g)$. Thurston, in his three page 1986 preprint ``A spine for
the \Teichmuller\ space of a closed surface'' \cite{Thurston:spine}, does construct a
model space for $\MCG(S_g)$. With some work, I can prove that the dimension of
Thurston's spine in genus~$g=2$ is indeed equal to $4g-5=8-5=3$. But I am unable to
prove that Thurston's spine in genus~$g\ge 3$ has dimension $4g-5$, and I think it may
be false. Ultimately we will depend on the Eilenberg-Ganea-Wall theorem
\cite{Brown:cohomology} to obtain an appropriate model space for
$\MCG(S_g)$, which is why we need to compute the $\vcd$.

\subsection{Mess subgroups}

Now we give Geoff Mess' proof of the lower bound: $\vcd(\MCG(S_g)) \ge 4g-5$.

The proof exhibits a subgroup $M_g \subgroup \MCG(S_g)$ which is a \Poincare\
duality group of dimension $4g-5$, in fact $M_g$ is the fundamental
group of a compact, aspherical $4g-5$ manifold. The group $M_g$ is called a \emph{Mess
subgroup} of $\MCG(S_g)$.

Mess subgroups are constructed by induction on genus.

\paragraph{Base case: Genus 2.} With $g=2$, we have $4g-5=3$, so we need a
3-dimensional subgroup of $\MCG(S_2)$. Take a curve family $\{c_1,c_2,c_3\} \subset
S_2$ consisting of three pairwise disjoint, pairwise nonisotopic curves. The Dehn
twists about $c_1,c_2,c_3$ generate a rank~3 free abelian group, and we are done.

Up to the action of $\MCG(S_2)$, there are two orbits of curve families
$\{c_1,c_2,c_3\}$, depending on whether or not some curve in the family separates. So,
there are two conjugacy classes of Mess subgroups in $\MCG(S_2)$.

\paragraph{Induction step:} Let $M_{g-1}$ be a Mess subgroup in $\MCG(S_{g-1})$, and
so $M_{g-1}$ is a \Poincare\ duality group of dimension $4(g-1)-5$. 

Consider the short exact sequence
$$1 \to \pi_1(S_{g-1}) \to \MCG(S_{g-1}^1) \to \MCG(S_{g-1}) \to 1
$$ 
Let $M'_{g-1} = $ preimage of $M_{g-1}$, so we get
$$1 \to \pi_1(S_{g-1}) \to M'_{g-1} \to M_{g-1} \to 1
$$
and it follows that $M'_{g-1}$ is \Poincare\ duality of dimension $4(g-1) - 5 + 2$. 

Let $S_{g,1}$ be the surface $S_g$ with a hole removed, and with one
boundary component. There is a central extension
$$1 \to \Z \to \MCG(S_{g,1}) \to \MCG(S_g^1) \to 1
$$
obtained by collapsing the hole to a punctuure; here, the group $\MCG(S_{g,1})$ is
defined as the group of homeomorphisms constant on the boundary, modulo isotopies
that are stationary on the boundary.

Let $M''_{g-1}$ be the preimage of $M'_{g-1}$ in $\MCG(S_{g,1})$, and we get
$$1 \to \Z \to M''_{g-1} \to M'_{g-1} \to 1
$$
from which it follows that $M''_{g-1}$ is \Poincare\ duality of dimension $4(g-1) - 5
+ 3$. 

Now attach a handle (a one-holed torus) to $S_{g,1}$ to get $S_{g+1}$, so
we get an embedding 
$$\MCG(S_{g,1}) \to \MCG(S_{g+1})
$$
Pick a simple closed curve $c$ contained in the handle and not isotopic to the
boundary. The Dehn twist $\tau_c$ commutes with $\MCG(S_{g,1})$, and in fact the
subgroup of $\MCG(S_{g+1})$ generated by $\tau_c$ and $\MCG(S_{g,1})$ is isomorphic to
the product $\MCG(S_{g,1}) \cross \tau_c$. We can therefore define
$$M_g = M''_{g-1} \cross \<\tau_c\>
$$
which is a \Poincare\ duality group of dimension $4(g-1) - 5 + 4 = 4g-5$ contained in
$\MCG(S_g)$.

This finishes Mess' proof that $\vcd(\MCG(S_g)) \ge 4g-5$.
\end{proof}

\paragraph{Remarks:} The construction of $M_g$ is completely determined by the isotopy
type of a certain filtration of $S_g$ by subsurfaces. There are only finitely many
such isotopy types up to the action of $\MCG$, and so there are only finitely many
conjugacy classes of Mess subgroups. In fact, a little thought shows that there are
exactly two conjugacy classes of Mess subgroups, distinguished by whether the
original curve system $\{c_1,c_2,c_3\}$ chosen in a genus~2 subsurface with one hole
contains a separating curve.

Letting $\Stab(c) \subgroup \MCG(S_g)$ be the stabilizer group of the closed curve $c$
picked in the last step, we have
$$M_g \subset \Stab(c)
$$
This fact will be significant later on.

\subsection{Model spaces and Mess cycles}

We will use a small trick: for the moment, we won't actually work with a model space
for $\MCG(S_g)$, instead we'll work with a model space for a finite index, torsion free
subgroup $\Gamma_g \subgroup \MCG(S_g)$. This is OK because the inclusion $\Gamma_g
\inject \MCG(S_g)$ is a quasi-isometry. Since $\vcd(\MCG(S_g)) = 4g-5$ it follows that
$\cd(\Gamma_g) = 4g-5$. We can therefore apply the Eilenberg-Ganea-Wall
theorem, to obtain a model space $B$ for $\Gamma_g$ of dimension $4g-5$. 

The reason for this trick is that for a group with torsion such as $\MCG(S_g)$, the
construction of a model space of dimension equal to the $\vcd$ is problematical.

Next we obtain a model space for $E$ and a fibration $\hyp^2 \to E \to B$ as follows.
Start with the canonical $\hyp^2$ bundle over \Teichmuller\ space $\Teich$, map $B$ to
$\Teich$ by a $\Gamma_g$-equivariant map, and pull the bundle back to get $E$. This
space $E$ is then a model space for the canonical $\pi_1(S_g)$ extension of
$\Gamma_g$, which is quasi-isometric to $\MCG(S_g^1)$. To prove that quasi-isometries
of $\MCG(S_g^1)$ coarse respect cosets of $\pi_1(S_g)$, it suffices to prove the same
for the $\pi_1(S_g)$ extension of $\Gamma_g$, and for this it suffices to prove that
quasi-isometries of $E$ coarsely respect the $\hyp^2$ fibers. Applying
Theorem~\ref{TheoremFibers}, it remains to construct top dimensional uniformly finite
cycles in $B$ which coarsely separate points.

Given a Mess subgroup $M \subgroup \MCG(S_g)$, the intersection $M' = M \intersect
\Gamma_g$ has finite index in $M$, and so $M'$ is still \Poincare\ duality of
dimension $4g-5$. The complex $B / M'$ is therefore a $K(M',1)$ space of dimension
$4g-5$. Since $M'$ is a \Poincare\ duality group, the homology $H_{4g-5}(M')$ is
infinite cyclic generated by the fundamental class $[M']$, and this class is
represented by a unique $4g-5$ cycle in $B/M'$. This cycle lifts to a $4g-5$
dimensional, uniformly finite cycle in $B$; call this a \emph{Mess cycle} in $B$.

We will prove that the Mess cycles in $B$ separate points.

\subsection{Passage to cosets of curve stabilizers}

We now pass from Mess cycles to left cosets of Mess subgroups to left cosets of curve
stabilizers, as follows. Although $\MCG$ does not act on $B$, it does quasi-act, which
is good enough. The quasi-action of $\MCG$ permutes the Mess cycles. There is a
bijection between Mess subgroups and Mess cycles: each Mess subgroup $M$ corresponds
to a unique Mess cycle $c$ such that $M$ that (coarsely) stabilizes $c$. If $M$
(coarsely) stabilizes $c$ and if $\Phi \in\MCG(S_g)$ then $\Phi M\Phi^\inv$ (coarsely)
stabilizes $\Phi(c)$.

Pick representatives $M_1,\ldots,M_k$ of the finitely many conjugacy classes of Mess
subgroups in $\MCG(S_g)$. It follows that, under the quasi-isometry
$B \to \MCG(S_g)$, Mess cycles correspond to left cosets in $\MCG(S_g)$ of
$M_1,\ldots,M_k$. So, it suffices to show that left cosets of $M_1,\ldots,M_k$
coarsely separate points in $\MCG(S_g)$.

Each Mess subgroup $M_i$ fixes some curve $c_i$, and so $M_i \subgroup
\Stab(c_i)$. Thus, each left coset of $M_i$ is contained in a left coset of
$\Stab(c_i)$. So, choosing curves $c_0,\ldots,c_n$ representing the orbits of simple
closed curves, it suffices to prove that the left cosets of the groups $\Stab(c_i)$
coarsely separate points in $\MCG$.

\subsection{New model space}

We now switch to a new model space $\Gamma$ for $\MCG(S_g)$, no longer contractible. We
will pass from left cosets of the groups $\Stab(c_i)$ to subsets of the new model
space~$\Gamma$.

$\Gamma$ is a graph whose vertices are pairs $(C,D)$ where each of $C,D$ is a pairwise
disjoint curve system, the systems $C,D$ jointly fill the surface, and each component
of $S-(C \union D)$ is a hexagon. This implies that $\MCG$ acts on the vertex set with
finitely many orbits. Since $\MCG$ is finitely generated, and since there are
finitely many orbits of vertices, it follows that we can attach edges in an
$\MCG$-equivariant way so that the graph $\Gamma$ is connected and has finitely
many orbits of edges. There's probably some nice scheme for attaching
edges, based on low intersection numbers, but it's not necessary. The
graph $\Gamma$ is now quasi-isometric to $\MCG$. Given a curve $c$, define $\Gamma_c$
to be the subgraph of $\Gamma$ spanned by vertices $(C,D)$ such that $c \in C \union
D$. 

\subsection{The subgraphs $\Gamma_c$ coarsely separate points}

We can now pass from left cosets of curve stabilizers to the sets $\Gamma_c$. Our
ultimate goal is to show that the system of subgraphs $\Gamma_c$, one for each curve
$c$, coarsely separates points in $\Gamma$.

Given vertices $(C,D)$ and $(C',D')$ which are very far from each other, we shall
pick a curve $c$ in $C \union D$ and show that $(C',D')$ is far from $\Gamma_c$.
This is enough, because $(C,D)$ is contained in $\Gamma_c$.

Since $(C,D)$ and $(C',D')$ are very far from each other, there exists $c \in C
\union D$ and $c' \in C' \union D'$ such that the intersection number $<c,c'>$
is very large. Proof: fixing $(C,D)$, if all such intersection numbers
$<c,c'>$ are uniformly small, then there is a uniform cardinality to the
number of possible $(C',D')$, so the distance from $(C,D)$ to $(C',D')$ is
uniformly bounded.

Consider now any curve system $(C_1,D_1)$ in $\Gamma_c$, meaning that $(C_1,D_1)$
contains $c$. The curve $c \in C_1 \union D_1$ has very large intersection
number with the curve $c' \in C' \union D'$. It follows that $(C_1,D_1)$ and
$(C',D')$ are far from each other. Proof: if $(C_1,D_1)$ and $(C',D')$ are
close, there is a uniform bound to the intersection number of a curve in
$C_1 \union D_1$ with a curve in $(C',D')$.

This completes the proof that quasi-isometries of $\MCG(S_g^1)$ coarsely
preserve fibers.

\section{Dynamical techniques: extensions of surface groups by pseudo-Anosov
homeomorphisms.}

In this section we give the proof of Theorem~\ref{TheoremMCGQIRigidity}, by proving
Strong QI-Rigidity of $\MCG(S_g^1)$ in the sense of Section~1. By applying
Theorem~\ref{TheoremCosets}, we are reduced to showing the following:
\begin{description}
\item[Fibered QI-rigidity] For all $K \ge 1$, $C \ge 0$, $R \ge 0$ there exists $A \ge
0$ such that if $\Phi$ is a $K,C$ quasi-isometry of $\MCG(S_g^1)$, and if $\Phi$ takes
each coset of $\pi_1 S_g$ to within a Hausdorff distance $R$ of some other coset of
$\pi_1 S_g$, then there exists $h \in \MCG(S_g^1)$ such that 
$$d_{\sup}(\Phi(f),h f) < A \quad\text{for all}\quad f \in \MCG(S_g^1)
$$
\end{description}

The methods of proof are very similar to the following result of Farb and myself:

\begin{theorem}[\cite{FarbMosher:sbf}]
\label{TheoremSchottky}
If $F$ is a Schottky subgroup of $\MCG(S_g)$ with extension
$$1 \to \pi_1(S_g) \to \Gamma_F \to F \to 1
$$
then the injection $\Gamma_F \to \QI(\Gamma_F)$ has finite index image. Moreover, if
$H$ is a group quasi-isometric to $\Gamma_F$ then there exists a homomorphism $H \to
\QI(\Gamma_F)$ with finite kernel and finite index image.
\end{theorem}

In that proof, using a tree as a model space for $F$, we proved that every
quasi-isometry of $\Gamma_F$ coarsely preserves fibers. Then we used pseudo-Anosov
dynamics (as we will here) to prove ``Fibered QI-rigidity'', from which
Theorem~\ref{TheoremSchottky} follows.

\subsection{\Teichmuller\ space and its canonical $\hyp^2$ bundle.} 
\label{SectionTeich}
We have already briefly mentioned these objects; here they are in more detail.

The \Teichmuller\ space $\Teich_g$ of $S_g$ is the space of hyperbolic structures on
$S_g$ modulo isotopy, or equivalently the space of conformal structures on $S_g$
modulo isotopy. We also need the \Teichmuller\ space $\T_g^1$ of
$S_g^1$, defined similarly using finite area complete hyperbolic structures on
$S_g^1$, or equivalently conformal structures with a removable singularity at the
puncture. It follows that each element of $\T_g^1$ can be expressed, up to isotopy, as
a pair $(\sigma,p)$ where $\sigma$ is a hyperbolic structure on $S_g$ representing an
element of $\T_g$, and $p$ is a point in the universal cover $\wt\sigma \approx
\hyp^2$. We therefore obtain a fiber bundle structure
$$\hyp^2 \to \T_g^1 \to \T_g
$$
on which the short exact sequence
$$1 \to \pi_1(S_g) \to \MCG(S_g^1) \to \MCG(S_g) \to 1
$$
acts. For each $x \in \T_g$ we use $\Sigma_x$ to denote the fiber of $\T_g^1$ over
$x$, and so $\Sigma_x$ is an isometric copy of $\hyp^2$. 

Note that the action of $\MCG(S_g^1)$ on $\T_g^1$ is \emph{not} cocompact, and so
we cannot regard $\T_g^1$ as a model space of $\MCG(S_g^1)$, and similarly for
$\MCG(S_g)$ acting on $\T_g$. There is, however, a cocompact
equivariant spine $Y_g \subset \T_g$
whose inverse image is a cocompact equivariant spine 
$Y_g^1 \subset \T_g^1$
and we have a fibration
$$\hyp^2 \to Y_g^1 \to Y_g
$$
on which the short exact sequence acts. By cocompactness, the
$Y$'s are model spaces for the $\MCG$'s.

The action of $\MCG(S_g^1)$ on $Y_g^1$ respects the $\hyp^2$ fibers, and the
stabilizer of each fiber is $\pi_1 S_g$. It follows that there is a quasi-isometry
$\MCG(S_g^1) \to Y_g^1$ taking cosets of $\pi_1 S_g$ to $\hyp^2$ fibers.

We can translate the result of Theorem~\ref{TheoremCosets} to the language of
$Y_g^1$. The translation says: every quasi-isometry $\Phi \from Y_g^1\to Y_g^1$
coarsely respects fibers, that is, there exists a constant $A\ge 0$ such that for each
$x\in Y_g$ there exists $x'\in Y_g$ such that $d_\Haus(\Phi(\Sigma_x),\Sigma_{x'}) \le
A$. Choosing an $x'$ for each $x$, we obtain an induced map $\phi \from Y_g \to Y_g$
with $\phi(x)=x'$, and $\phi$ is a quasi-isometry.

This reduces the proof of Fibered QI-rigidity for $\MCG(S_g^1)$ to the analogous
statement for quasi-isometries of $Y_g^1$: with $\Phi$ as above, we must find $h \in
\MCG(S_g^1)$ so that the actions of $\Phi$ and $h$ on $S_g^1$ agree to within bounded
distance.

\paragraph{$\hyp^2$ bundles over lines.} Consider bi-infinite, proper paths $\ell\from
\reals \to \T_g$, with image often in $Y_g$. The path $\ell$ is always piecewise smooth
and Lipschitz. Let 
$$\Sigma_\ell = \pi^\inv(\ell) \subset \T_g^1
$$
so we have an $\hyp^2$ bundle over the line $\ell$:
$$\hyp^2 \to \Sigma_\ell \to \ell
$$ 
There is a reasonably natural metric on $\Sigma_\ell$, obtained by combining the
$\hyp^2$ metric on fibers with the metric on the $\reals$ factor. There are some
choices, but the metric is natural up to quasi-isometry. If $\ell$ is piecewise
geodesic then for the metric on $\Sigma_\ell$ we can take the pullback of the metric
on $\T_g^1$.

Given a quasi-isometry $\Phi \from Y_g^1 \to Y_g^1$ with induced quasi-isometry $\phi
\from Y_g \to Y_g$, for any bi-infinite proper path $\ell \from \reals \to \T_g$ the
map $\Phi$ restricts to a fiber respecting quasi-isometry
$$\Sigma_\ell \to \Sigma_{\phi(\ell)}
$$
Also, if $\ell,\ell' \from \reals \to Y_g$ are ``fellow travellers'', then there is an
induced $\pi_1 S$-equivariant map $\Sigma_\ell \to \Sigma_{\ell'}$ which is a
quasi-isometry.

Example: Suppose that $\ell$ is the axis of a pseudo-Anosov diffeomorphism,
or more generally, that $\ell$ fellow travels such an axis. Thurston's
hyperbolization theorem for fibered 3-manifolds implies that $\Sigma_\ell$ is
quasi-isometric to $\hyp^3$, and so $\Sigma_\ell$ is a Gromov hyperbolic metric space.

\paragraph{\Teichmuller\ geodesics and their singular \solv\ spaces.} A \emph{quadratic
differential} on $S_g$ is a transverse pair of measured foliations:
$$q=(\fol_u,\fol_s)
$$
Each quadratic differential $q$ determines a singular Euclidean metric, which
determines conformal structure with removable singularities, which determines
a point $\sigma(q) \in \T_g$. For each $t \in \reals$ define 
$$q_t = (e^{-t} \fol_u, e_t \fol_s)
$$
The path 
$$\gamma_q = \{t \mapsto \sigma(q_t) \suchthat t \in \reals\}
$$ 
in $\T_g$ is the \emph{\Teichmuller\ geodesic} corresponding to $q$.

Let $\Sigma^\solv_q$ denote the hyperbolic plane bundle
$\Sigma_{\gamma_q}$ with the \emph{singular \solv\ metric}, defined by
$$e^{-2t} \, d\fol_u^2 + e^{2t} \, d\fol_s^2 + dt^2
$$

Fact: if $\gamma_q$ is cobounded in $\T_g$, meaning that it is contained in the
$\MCG(S_g)$ orbit of some bounded subset $B$ of $\T_g$, then the the identity map
$\Sigma_{\gamma_q}\to \Sigma^\solv_q$ is a quasi-isometry between the ``natural''
metric and the singular \solv\ metric, with quasi-isometry constants depending only on
$B$.

\paragraph{Pseudo-Anosov homeomorphisms, their axes, and their hidden
symmetries.} \hfill\break A homeomorphism $f \from S_g \to S_g$ is \emph{pseudo-Anosov}
if there exists a quadratic differential $q_f=(\fol_u,\fol_s)$ and $\lambda>1$ such
that 
$$f(\fol_u,\fol_s) = (\lambda^\inv \fol_u, \lambda \fol_s)
$$
It follows that the path $\gamma_{q_f}$ is invariant under $f$ in $\Teich$, and in
fact $\gamma_{q_f}$ is the set of points $\sigma \in \Teich$ at which $d(f
\sigma,\sigma)$ is minimized; we call $\gamma_{q_f}$ the \emph{axis} of $f$ in
$\Teich$. Conversely, a mapping class $\Phi\in\MCG(S_g)$ is represented by a
pseudo-Anosov homeomorphism only if $d(\Phi\sigma,\sigma)$ has a positive
minimum in $\T$. These facts were proved by Bers \cite{Bers:ThurstonTheorem}. 

Let $\Sigma_f$ denote $\Sigma_{q_f}$, with a superscript ``\solv'' added to the
notation is we wish to denote the singular \solv\ metric. The group $J_f =
\pi_1(S_g)\semidirect_f\Z$ acts by isometries on $\Sigma^\solv_f$, that is, 
$$J_f \subgroup I_f = \Isom(\Sigma^\solv_f))
$$
It is possible that $J_f$ is properly contained in $I_f$. We can think
of the elements of $I_f - J_f$ as ``hidden symmetries'' of $f$. One possibility for
hidden symmetries occurs when $f$ is a proper power. Another possibility occurs when
$f$ or some power of $f$ is conjugate to its own inverse. In general, $I_f$ is a
virtually cyclic group, and like all such groups there is an epimorphism $I_f \to C$
whose image $C$ is either infinite cyclic or infinite dihedral, and whose kernel is
finite; a nontrivial kernel provides a further source of hidden symmetries
of $f$. 

We shall need an alternate description of $I_f$. There exists a maximal index orbifold
subcover $S_g \to O_f$ such that $f$ descends to a pseudo-Anosov homeomorphism of the
orbifold~$O_f$ which we shall denote $f'$. Let $\VN_{f'}$ be the virtual normalizer of
$f'$ in $\MCG(\O_f)$, consisting of all $g \in \MCG(\O_f)$ such that $g^\inv \< f' \> g
\intersect \<f'\>$ has finite index in each of the infinite cyclic subgroups $g^\inv
\<f'\> g$ and $\<f'\>$. 

\begin{fact} There is a natural extension
$$1 \to \pi_1(\O_f) \to I_f \to \VN_{f'} \to 1
$$
\end{fact}

Since the virtual normalizer of a pseudo-Anosov homeomorphism is virtually cyclic, it
follows that $I_f$ contains $J_f$ with finite index.

A direct construction can be used to show:

\begin{fact}
There exist pseudo-Anosov homeomorphisms with no hidden symmetries, that is,
so that $I_f = J_f$.
\end{fact}

For example, consider the fact the dimension of the measured foliation space of $S_g$
is $6g-6$, and the transition matrix of a train track representative of every
pseudo-Anosov homeomorphism $f$ is an $n \cross n$ matrix with $n \le 6g-6$; it follows
that the algebraic degree of the expansion factor $\lambda(f)$ is at most $6g-6$. One
can construct primitive pseudo-Anosov homeomorphisms $f$ of $S_g$ such that
$\lambda(f)$ has maximal algebraic degree $6g-6$. For such an $f$, the kernel $K$ of
the epimorphism $I_f \to C$ must be trivial, for if $K$ is nontrivial then $f$ or some
power of $f$ commutes with $K$ and so descends through an orbifold covering map $S_g\to
\O=S_g/K$ where the measured foliation space for $\O$ has strictly smaller dimension
than for $S_g$; it follows that the degree of $\lambda(f)$ is strictly smaller than
$6g-6$. 

This still leaves open the possibility that $I_f$ is infinite dihedral, implying that
$f$ is conjugate to its own inverse, but a random example will fail to have this
property.

\subsection{Proof of fibered QI-rigidity} 

Consider a quasi-isometry $\Phi \from Y_g^1 \to Y_g^1$. As noted in
Section~\ref{SectionTeich}, $\Phi$ coarsely respects the fibers of the fibration
$Y_g^1 \to Y_g$, and so $\Phi$ induces a quasi-isometry $\phi \from Y_g \to Y_g$.
For any bi-infinite, proper path $\ell \from \reals \to Y_g$, the quasi-isometry
$\Phi$ induces a coarse fiber respecting quasi-isometry from $\Sigma_\ell$ to
$\Sigma_{\phi(\ell)}$. We shall apply various fiber respecting quasi-isometry
invariants to the metric spaces $\Sigma_\ell$. 

For example, we say that $\ell$ and its $\hyp^2$ bundle $\Sigma_\ell$ are
\emph{hyperbolic} if $\Sigma_\ell$ is a $\delta$-hyperbolic metric space for some
$\delta \ge 0$. Since hyperbolicity is a quasi-isometry invariant, we obtain:

\begin{fact}[Hyperbolic spaces are preserved]
\label{FactHyp}
Every quasi-isometry $\Phi$ of $\MCG(S_g^1)$ coarsely respects the hyperbolic
spaces $\Sigma_\ell$, and the induced quasi-isometry $\phi$ of $\MCG(S_g)$ coarsely
respects the hyperbolic paths $\ell$.
\end{fact}

Much more surprising is that a quasi-isometry $\Phi$ coarsely respects the
\emph{periodic} hyperbolic spaces. To be precise, a \Teichmuller\ geodesic $\gamma
\subset \T_g$, or its $\hyp^2$ bundle $\Sigma_\gamma$, is \emph{periodic} if $\gamma$ is
the axis of some pseudo-Anosov mapping class in $\MCG(S_g)$. Equivalently, by
Thurston's hyperbolization theorem, $\Sigma_\gamma$ is the universal cover of a fibered
hyperbolic 3-manifold. A bi-infinite path $\ell \from \reals \to Y_g$ is
\emph{coarsely periodic} if there exists an axis $\gamma$ in $\Teich$ such that $\ell$
and $\gamma$ are fellow travellers, meaning that $d(\ell(h(t)),\gamma(t))$ is uniformly
bounded where $h \from \reals \to \reals$ is some quasi-isometry of the real line.

Here is the heart of the matter:

\begin{theorem}
\label{TheoremPeriodic}
Every quasi-isometry of $\MCG(S_g^1)$ coarsely respects the periodic
hyperbolic 3-manifolds $\Sigma_\gamma$. To be precise, given: a quasi-isometry $\Phi
\from Y_g^1 \to Y_g^1$ inducing a quasi-isometry $\phi \from Y_g \to Y_g$, given
a periodic axis $\gamma \subset \T_g$, and given a coarsely periodic path $\ell$ in
$Y_g$ that fellow travels $\gamma$, there exists a periodic axis $\gamma' \subset \T_g$
such that $\phi \composed \ell$ fellow travels $\gamma'$. Moreover, $\Phi \from
\Sigma_\ell\to \Sigma_{\phi(\ell)}$ is a bounded  distance from an isometry
$H \from \Sigma^\solv_\gamma \to \Sigma^\solv_{\gamma'}$. 
\end{theorem}

The meaning of the final sentence of this theorem is that there is a commutative
diagram of fiber respecting quasi-isometries
$$\xymatrix{
\Sigma_\ell \ar[r]^{\Phi} \ar[d] & \Sigma_{\phi(\ell)} \ar[d] \\
\Sigma^\solv_\gamma \ar[r]_H & \Sigma^\solv_{\gamma'}
}
$$
where the vertical arrows represent maps that move points a uniformly bounded distance
in $\T_g^1$.

In this theorem, the bounds in the conclusion depend only on the quasi-isometry
constants of $\Phi$ and on the fellow traveller constants for $\ell$ and $\gamma$. 

Before proving this theorem we first apply it to:

\begin{proof}[Proof of QI-rigidity of $\MCG(S_g^1)$] \quad

Fix a quasi-isometry $\Phi\from Y_g^1 \to Y_g^1$. 

Consider a pseudo-Anosov homeomorphism $f \from S_g \to S_g$ without hidden symmetries:
the group $J_f = \pi_1(S_g) \semidirect_f \Z$ is the entire isometry group of the
singular \solv\ manifold $\Sigma^\solv_f$. It follows that each conjugate $gfg^\inv$
has no hidden symmetries. By Theorem~\ref{TheoremPeriodic}, there is a pseudo-Anosov
$f' \from S_g \to S_g$ and an isometry $H_{f} \from \Sigma^\solv_f
\to \Sigma^\solv_{f'}$ such that $\Phi$ takes $\Sigma_f$ to $\Sigma_{f'}$ by a map
which is a bounded distance from $H_f$. The isometry $H_f$ conjugates
$I_f =\Isom(\Sigma^\solv_f)$ to $I_{f'} = \Isom(\Sigma^\solv_{f'})$, implying that
$H_f$ conjugates $\pi_1(\O_f)$ to $\pi_1(\O_{f'})$. However, since $f$ has no
hidden symmetries, $\O_f = S_g$, and so $\pi_1(\O_{f'})$ must also equal $S_g$;
replacing $f'$ if necessary by some root, it follows that $f'$ has no hidden
symmetries. The isometry $H_f$ therefore agrees with the action of some automorphism
of $\pi_1(S_g)$, which is identified with a mapping class $h_f \in \MCG(S_g^1)$, and
conjugation by $h_f$ takes $f$ to $f'$.

Next we must show that $h_f$ is independent of $f$. For this we use the well known
fact that $\MCG(S_g^1)$ acts faithfully on the circle at infinity $S^1_\infinity$ of
$\hyp^2 = \wt S_g$; this fact is the basis of Nielsen theory. 

The quasi-isometry $\Phi$ also acts on $S^1_\infinity$. To see why, fix a fiber
$\Sigma_x$ and identify this fiber isometrically with $\hyp^2$, so that the boundary
of $\Sigma_x$ is identified with $S^1_\infinity$. Since $\Phi(\Sigma_x)$ is
(uniformly) coarsely equivalent to some fiber $\Sigma_{x'}$, and since $\Sigma_{x'}$
is (nonuniformly) coarsely equivalent to $\Sigma_x$, the action of $\Phi$ induces a
self quasi-isometry of $\Sigma_x$, thereby inducing a homeomorphism of $S^1_\infinity$.

By construction of $h_f$, the actions of $\Phi$ and of $h_f$ on $S^1_\infinity$ agree,
that is to say, the action of $h_f$ on $S^1_\infinity$ is independent of $f$. By
faithfulness of the action of $\MCG(S_g^1)$ on $S^1_\infinity$, it follows that $h=h_f
\in \MCG(S_g^1)$ is independent of $f$. We have thus constructed the desired mapping
class $h \in \MCG(S_g^1)$, and from the construction it is evident that $\Phi \from
Y_g^1 \to Y_g^1$ is within bounded distance of the action of $h$.

\end{proof}

\subsection{Proof of Theorem \ref{TheoremPeriodic}}

This theorem reduces quickly to results from \cite{Mosher:StableQuasigeodesics}
(see also \cite{Bowditch:stacks}) and from \cite{FarbMosher:sbf}, for each of which we
will sketch a proof in broad strokes.

\paragraph{Step 1: Hyperbolic lines in $\T_g$.} First we need the following theorem,
proved independently by Bowditch and by myself:

\begin{theorem}[\cite{Bowditch:stacks}; \cite{Mosher:StableQuasigeodesics}]
\label{TheoremCoboundedGeodesics}
A line $\ell \from \reals \to \T_g$ is hyperbolic if and only if there exists
a cobounded \Teichmuller\ geodesic $\gamma \from \reals \to \T_g$ that fellow travels
$\ell$.
\end{theorem}

For example, a pseudo-Anosov axis is a cobounded \Teichmuller\ geodesic,
forming a countable family. In toto, there are uncountably many cobounded
\Teichmuller\ geodesics, making Theorem~\ref{TheoremPeriodic} all the more surprising.

\begin{proof}[A one minute proof] (Every theorem should have a one minute proof,
a five minute proof, a twenty minute proof\ldots) 

We must construct a quadratic differential $q = (\fol_s,\fol_u)$. 

``Hyperbolicity'' means ``exponential divergence of
geodesics'' \cite{Cannon:TheoryHyp}. Ordinarily this applies to geodesic rays passing
transversely through spheres, but it also applies to geodesics in $\Sigma_\ell$
passing transversely through fibers the $\Sigma_t = \ell^\inv(t)$
\cite{FarbMosher:quasiconvex}. It follows that every geodesic contained in a fiber
$\Sigma_t$ is stretched exponentially in either the forward or backward direction, as
$t \to +\infinity$ or as $t \to -\infinity$. Some geodesics are stretched exponentially
in both directions; indeed, this is true of a random geodesic in a fiber. Certain
geodesics contained in the fibers $\Sigma_t$ are stretched exponentially as
$t\to\infinity$, but not as $t \to -\infinity$; these geodesics form the unstable
foliation
$\fol_u$. Certain other geodesics in $\Sigma_t$ are stretched exponentially as $t
\to -\infinity$ but not as $t \to +\infinity$, and these form the stable foliation
$\fol_s$. 

Taking $q=(\fol_u,\fol_s)$, a compactness argument shows that $\ell$ fellow travels
the \Teichmuller\ geodesic $\gamma_q$.
\end{proof}

\paragraph{Step 2: Periodic hyperbolic lines.} The key fact, from which
Theorem~\ref{TheoremPeriodic} quickly follows, is: 

\begin{theorem}[\cite{FarbMosher:sbf}] 
\label{TheoremSolvQI}
Let $\gamma,\gamma'$ be cobounded geodesics in
$\T_g$, and suppose that $\gamma$ is periodic. If there exists a fiber respecting
quasi-isometry $\Phi\from\Sigma^\solv_\gamma\to\Sigma^\solv_{\gamma'}$ then $\gamma'$
is periodic and $\Phi$ is a bounded distance from an isometry.
\end{theorem}

Before sketching the proof, we apply it to:

\begin{proof}[Proof of Theorem \ref{TheoremPeriodic}]
We replace the map $\MCG(S_g^1) \to \MCG(S_g)$, fibered by cosets of $\pi_1 S_g$, with
the map $Y_g^1 \to Y_g$, fibered by copies of $\hyp^2$. Let $\Phi \from Y_g^1 \to
Y_g^1$ be a quasi-isometry, inducing a quasi-isometry $\phi \from Y_g \to Y_g$. Let
$\gamma \subset \T_g$ be a periodic axis fellow travelling a coarsely periodic path
$\ell$. We obtain a fiber respecting quasi-isometry $\Sigma_\ell \to
\Sigma^\solv_\gamma$. The quotient $\Sigma^\solv_\gamma$ modulo its isometry group is
a hyperbolic \nb{3}-orbifold, by Thurston's hyperbolization theorem, and so
$\Sigma^\solv_\gamma$ is quasi-isometric to $\hyp^3$. It follows that $\Sigma_\ell$ is
a hyperbolic metric space, that is, $\ell$ is a hyperbolic line. By
Fact~\ref{FactHyp}, $\phi(\ell)$ is a hyperbolic line and $\Sigma_{\phi(\ell)}$
is a hyperbolic metric space. Applying Theorem~\ref{TheoremCoboundedGeodesics}, there
is a cobounded geodesic $\gamma'$ in $\Teich$ that fellow travels $\phi(\ell)$. By
combining the fiber respecting quasi-isometries $\Sigma^\solv_\gamma \to \Sigma_\ell
\to \Sigma_{\Phi(\ell)} \to \Sigma^\solv_{\gamma'}$, we obtain a fiber respecting
quasi-isometry $\Sigma^\solv_\gamma \to \Sigma^\solv_{\gamma'}$. Applying
Theorem~\ref{TheoremSolvQI}, $\gamma'$ is periodic the latter quasi-isometry is a bounded
distance from an isometry. 
\end{proof}

\paragraph{Step 3: Proof of Theorem \ref{TheoremSolvQI}.}  Farb and I gave a proof
that uses Thurston's hyperbolization theorem for fibered 3-manifolds
\cite{Otal:fibered} together with Rich Schwartz' geodesic pattern rigidity theorem
\cite{Schwartz:Symmetric}. It would be extremely nice to have a proof which uses only
pseudo-Anosov dynamics, but I still don't know how to do this. Here is a broad sketch
of the proof of Farb and myself taken from \cite{FarbMosher:sbf}.

Given a cobounded \Teichmuller\ geodesic $\gamma$, define $\QI_f(\Sigma^\solv_\gamma)$
to be the group of ``fiber respecting quasi-isometries'' of $\Sigma^\solv_\gamma$. We
have an injection $\Isom(\Sigma^\solv_\gamma) = I_\gamma \inject
\QI_f(\Sigma^\solv_\gamma)$, and the question arises whether there is anything else in
$\QI_f(\Sigma^\solv_\gamma)$.

First we prove, when $\gamma$ is a periodic geodesic, that the injection
$I_\gamma \to \QI_f(\Sigma^\solv_\gamma)$ is an isomorphism,
that is, every self quasi-isometry of $\Sigma^\solv_\gamma$ that coarsely respects
fibers is a bounded distance from an isometry. 

The singular lines of $\Sigma^\solv_\gamma$ form a collection of singular \solv\
geodesics intersecting the fibers at right angles; let $\Omega$ denote this collection
of geodesics. If $\Phi$ is a fiber respecting quasi-isometry of $\Sigma^\solv_\gamma$,
then $\Phi$ coarsely respects leaves of $f^s$ and $f^u$ as noted earlier, and in fact
$\Phi$ coarsely respects the suspensions of these leaves. It follows that $\Phi$
coarsely respects $\Omega$, because the singular lines in $\Omega$ are precisely the
sets which, coarsely, are intersections of three or more suspensions of leaves of
$f^s$ (or of $f^u$) whose pairwise intersections are unbounded. We may then move
values of $\Phi$ by a bounded amount so that $\Phi$ is a homeomorphism that strictly
respects leaves of $f^s$, leaves of $f^u$, and $\Omega$.

By Thurston's hyperbolization theorem, there is an $I_\gamma$-equivariant
$\hyp^3$ metric on $\Sigma^\solv_\gamma$. The lines in $\Omega$ can be straightened
to hyperbolic geodesics, which are evidently invariant under $I_\gamma$. This
is exactly the setup of Schwartz' theorem, whose conclusion is that the
group of quasi-isometries of $\hyp^3$ that coarsely respects $\Omega$ contains
$I_\gamma$ with finite index, that is, $\QI_f(\Sigma^\solv_\gamma)$ contains
$I_\gamma$ with finite index. But then an easy argument shows that $I_\gamma$ must
actually be all of $\QI_f(\Sigma^\solv_\gamma)$.

Using the isomorphism $I_\gamma \to \QI_f(\Sigma^\solv_\gamma)$ and the
coarse fiber respecting quasi-isometry $\Phi \from \Sigma^\solv_\gamma \to
\Sigma^\solv_{\gamma'}$, we now show that $\Phi$ is a bounded distance from an
isometry. As above, we may first move values of $\Phi$ by a bounded amount so that
$\Phi$ is a homeomorphism that respects the fibers and the stable and unstable
foliations. Let $\Sigma_x$ denote a fiber of $\Sigma^\solv_\gamma$ whose image under
$\Phi$ is a fiber $\Sigma_{x'}$ of $\Sigma^\solv_{\gamma'}$. Let $I_x$ be the subgroup
of $\Isom(\Sigma_x)$ that preserves the stable foliation and the unstable foliation,
and similarly for $I_{x'}$. The conjugate action $\Phi^\inv \composed I_{x'} \composed
\Phi$ is an action on $\Sigma^\solv_\gamma$ by quasi-isometries that preserves each
fiber. By the computation of $\QI_f(\Sigma^\solv_\gamma)$ just given, we obtain
$\Phi^\inv \composed I_{x'} \composed \Phi \subset I_x$, and in particular $\Phi^\inv
\composed I_{x'} \composed \Phi$ preserves the invariant measures on the stable and
unstable foliations of $\Sigma_x$. Conjugating back now to $\Sigma_{x'}$, it follows
that $I_{x'}$ preserves two sets of invariant measures on the stable and unstable
foliations of $\Sigma_{x'}$: the ones coming from the singular \solv\ structure on
$\Sigma^\solv_{\gamma'}$, and the ones pushed forward via $\Phi$. But the stable and
unstable foliations on $\Sigma_{x'} / I_{x'}$ are uniquely ergodic: this follows from
a theorem of Masur \cite{Masur:UniquelyErgodic}, which says that the stable and
unstable foliations associated to a cobounded geodesic in \Teichmuller\ space are
uniquely ergodic. Up to rescaling, therefore, the map that $\Phi$ induces from
$\Sigma_x$ to $\Sigma_{x'}$ is an isometry; and the rescaling may be ignored by moving
the point $x'$ up or down. But this immediately implies that $\Sigma^\solv_\gamma$ and
$\Sigma^\solv_{\gamma'}$ are isometric, by an isometry that agrees with $\Phi$ on
$\Sigma_x \to \Sigma_{x'}$ and that moves other fibers up or down by uniformly bounded
adjustements.


\providecommand{\bysame}{\leavevmode\hbox to3em{\hrulefill}\thinspace}

\noindent
Lee Mosher:\\
Department of Mathematics, Rutgers University, Newark\\
Newark, NJ 07102\\
E-mail: mosher@andromeda.rutgers.edu

\end{document}